\documentclass[12pt,headings=normal]{scrartcl}

\usepackage[utf8]{inputenc}
\usepackage[ngerman]{babel}
\usepackage{amsmath,amssymb}
\usepackage[thmmarks]{ntheorem}
\usepackage{verbatim}
\usepackage{graphicx}
\usepackage[all]{xy}

\newcommand{\Ass}{\operatorname{Ass}}
\newcommand{\Ann}{\operatorname{Ann}}

\newcommand{\Hom}{\operatorname{Hom}}

\theorembodyfont{\itshape}
\theoremstyle{plain}
\newtheorem{Satz}{Satz}[section]
\newtheorem{Folgerung}[Satz]{Folgerung}
\newtheorem{Lemma}[Satz]{Lemma}
\newtheorem{Proposition}[Satz]{Proposition}
\theorembodyfont{\upshape}
\newtheorem{Bemerkung}[Satz]{Bemerkung}

\theoremheaderfont{\scshape}
\theorembodyfont{\upshape}
\theoremstyle{nonumberplain}
\theoremseparator{}
\theoremsymbol{$\Box$}
\newtheorem{Beweis}{Beweis}

\title{\Large \"Uber die von einem Ideal $I \subset R$ erzeugten $R$-Moduln III}
\author{\large Helmut Z\"oschinger\\
  \large Mathematisches Institut der Universit\"at M\"unchen\\
  \large Theresienstr. 39, D-80333 M\"unchen\\
  \large E-mail: zoeschinger$@$mathematik.uni-muenchen.de
}
\date{}

\newcounter{abccount}
\newenvironment{abc}{%
\begin{list}{(\alph{abccount})}{%
\usecounter{abccount}%
\setlength{\partopsep}{0pt}%
\setlength{\topsep}{1ex}%
\setlength{\itemsep}{0pt}%
}%
}{\end{list}}

\newcounter{myenumcount}
\newenvironment{myenum}{%
\begin{list}{(\arabic{myenumcount})}{%
\usecounter{myenumcount}%
\setlength{\partopsep}{0pt}%
\setlength{\topsep}{1ex}%
\setlength{\itemsep}{0pt}%
}%
}{\end{list}}

\newcounter{iiicount}
\newenvironment{iii}{%
\begin{list}{(\roman{iiicount})}{%
\usecounter{iiicount}%
\setlength{\labelwidth}{3em}%
\setlength{\partopsep}{0pt}%
\setlength{\topsep}{1ex}%
\setlength{\itemsep}{0pt}%
}%
}{\end{list}}

\begin{document}
\maketitle

\centerline{\textbf{Abstract}}
\begin{abstract}
  \noindent
  Let $(R, \mathfrak m)$ be a commutative noetherian local ring and $I$ an
  ideal of $R$. For every $R$-module $M$, $\gamma_I(M) = \sum\{
  \operatorname{Bi} f \,|\, f \in \operatorname{Hom}_R(I,M)\}$ is called the
  \emph{trace} of $I$ in $M$. It is easy to see that
  $\operatorname{Ext}_R^1(R/I,M) = 0$ always implies $IM = \gamma_I(M)$. If
  the second condition holds for all ideals $I$ of $R$, we say that $M$ is
  \emph{excellent}. In part~1, we show a number of conditions for these
  modules, which are well-known for injective modules. In the second part,
  we examine the special case $M = R$. In particular, we show that for every
  prime ideal $\mathfrak{p}$ the equality $\mathfrak{p} =
  \gamma_{\mathfrak{p}}(R)$ holds iff $R_{\mathfrak{p}}$ is \emph{not} a
  discrete valuation ring. From the results by Matlis (1973) about
  1-dimensional local CM-rings and with the help of the first neighborhood
  ring $\Lambda$, it follows immediately that $\gamma_{\mathfrak{m}^n} (R) =
  \Lambda^{-1}$ for almost all $n \geq 1$. In the third part, we examine the
  dual construction $\kappa_I(M) = \bigcap \{ \operatorname{Ke} f \,|\, f\in
  \operatorname{Hom}_R(M,I^\circ) \}$ and reduce the main results about
  $\operatorname{Tor}_1^R(M, R/I) = 0$ and $\kappa_I(M) = M[I]$ to part~1 by
  considering the Matlis dual $M^\circ = \operatorname{Hom}_R(M, E)$ and the
  equalities $\gamma_I(M^\circ) = \operatorname{Ann}_{M^\circ}(\kappa_I(M))$,
  $\kappa_I(M^\circ) = \operatorname{Ann}_{M^\circ}(\gamma_I(M))$.
\end{abstract}

\medskip


\bigskip


\section{Der Untermodul $\gamma_I(M)$}

Stets sei $(R, \mathfrak{m})$ ein kommutativer, noetherscher, lokaler Ring,
$M$ ein $R$-Modul und $I$ ein Ideal von $R$. Wie in \cite{018} bezeichnen
wir mit $\gamma_I(M)$ den gr\"o{\ss}ten $I$-generierten Untermodul von $M$, d.\,h.
$\gamma_I(M) = \sum \{\operatorname{Bi} f \,|\, f \in
\operatorname{Hom}_R(I,M) \}$, und damit gilt stets
\begin{equation*}
  IM \subset \gamma_I(M) \subset M[\Ann_R(I)].
\end{equation*}

Falls $M$ injektiv ist, stimmen alle drei Untermoduln \"uberein (\cite{018}
Folgerung 3.2), falls $I$ zyklisch ist (es gen\"ugt $I \twoheadrightarrow
R/\Ann_R(I)$, siehe \cite{018}~p.~5), gilt $\gamma_I(M) = M[\Ann_R(I)]$. Uns
interessiert im Folgenden die erste Ungleichung, und wir sagen $M$ sei
\emph{$I$-ausgezeichnet}, wenn $IM = \gamma_I(M)$ gilt.

Der Extremfall $\gamma_I(M) = 0$ ist einfach:
\begin{Lemma}\label{1.1}
  F\"ur einen $R$-Modul $M$ und ein Ideal $I \subset R$ sind \"aquivalent:
  \begin{iii}
    \item $\gamma_I(M) = 0$.
    \item $\operatorname{Hom}_R(I,M) = 0$.
    \item $I_{\mathfrak{p}} = 0$ f\"ur alle $\mathfrak{p} \in \Ass(M)$.
  \end{iii}
\end{Lemma}

\begin{Beweis}
  Nach (\cite{018} Satz~3.1,~c) ist $\Ass(\gamma_I(M)) = \Ass(M) \cap
  \operatorname{Supp}(I)$, nach (\cite{002} chap. IV, §\,1, Prop.~10) aber
  die rechte Seite gleich $\Ass(\operatorname{Hom}_R(I,M))$, und damit folgt
  alles. 
\end{Beweis}

\begin{Proposition}\label{1.2}
  F\"ur einen $R$-Modul $M$ und ein Ideal $I \subset R$ sind \"aquivalent:
  \begin{iii}
  \item
    $M$ ist $I$-ausgezeichnet, d.\,h. $IM = \gamma_I(M)$.
  \item F\"ur jeden Homomorphismus $f \colon I \to M$ gilt $\operatorname{Bi}
    f \subset IM$, d.\,h.
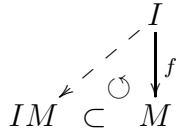
    \begin{equation*}
      \xymatrix{
        \ar@{}|(.7){\textstyle{\circlearrowleft}}[dr]&I\ar@{-->}[dl]\ar[d]^f\\
        IM\ar@{}|{\textstyle\subset}[r]&M
      }
    \end{equation*}
  \item $\operatorname{Hom}_R(I,M) \xrightarrow{\ \nu_{\bullet}\ }
    \operatorname{Hom}_R(I, M/IM)$ ist die Nullabbildung.
  \end{iii}
  War das Ideal $I$ zyklisch, ist das weiter \"aquivalent mit
  \begin{iii}
  \setcounter{iiicount}{3}  
  \item $IM = M[\Ann_R(I)]$.
  \end{iii}
\end{Proposition}

\begin{Beweis}
  (ii $\leftrightarrow$ iii) ist klar, ebenso (i $\to$ ii), denn f\"ur jedes
  $f \in \operatorname{Hom}_R(I,M)$ ist $\operatorname{Bi} f$ $I$-generiert,
  also nach Voraussetzung enthalten in $IM$.\\
  Bei (ii $\to$ i) ist nur "`$\supset$"' zu zeigen, und mit $g = \langle
  g_{\lambda} \rangle \colon I^{(\Lambda)} \twoheadrightarrow \gamma_I(M)$
  ist nach Voraussetzung $\operatorname{Bi} g_{\lambda} \subset IM$ f\"ur alle
  $\lambda \in \Lambda$, also auch $\sum \operatorname{Bi} g_{\lambda} =
  \operatorname{Bi} g = \gamma_I(M)$
  enthalten in $IM$.\\
  (iv $\to$ i) gilt nach der Einleitung stets, wenn aber $I$ zyklisch ist,
  sogar $\gamma_I(M) = M[\Ann_R(I)]$, also auch (i $\to$ iv).
\end{Beweis}

\begin{Bemerkung}
  Im Allgemeinen ist die Bedingung (iv) st\"arker als (i): Genau dann ist $R$
  $\mathfrak{m}$-ausgezeichnet, wenn $R$ \emph{kein} diskreter
  Bewertungsring ist. Falls also $R$ ein Integrit\"atsring mit $\dim(R) \geq
  2$ ist, wird $R[\Ann_R(\mathfrak{m})]$ echt gr\"o{\ss}er als
  $\mathfrak{m}$.\hfill$\Box$ 
\end{Bemerkung}

\begin{Folgerung}\label{1.4}
  Die Familie aller $I$-ausgezeichneten $R$-Moduln ist gegen\"uber direkten
  Produkten, reinen Untermoduln und direkten Summen abgeschlossen.
\end{Folgerung}

\begin{Folgerung}\label{1.5}
  Ist $(I_{\lambda} \,|\, \lambda \in \Lambda)$ eine Familie von Idealen und
  $M$ $I_{\lambda}$-ausgezeichnet f\"ur alle $\lambda \in \Lambda$, so ist $M$
  auch $(\sum I_{\lambda})$-ausgezeichnet.
\end{Folgerung}

\begin{Proposition}\label{1.6}
  F\"ur einen $R$-Modul $M$ und ein Ideal $I \subset R$ sind \"aquivalent:
  \begin{iii}
    \item $\operatorname{Ext}_R^1(R/I, M) = 0$.
    \item $M$ ist $I$-ausgezeichnet und die kanonische Abbildung $\sigma
      \colon M \to \operatorname{Hom}_R(I,IM)$ surjektiv.
  \end{iii}
\end{Proposition}

\begin{Beweis}
  $\sigma$ ist eine Verallgemeinerung der wohlbekannten Abbildung $R \to
  \operatorname{End}_R(I)$, $r \mapsto \tilde{r}$ und ist definiert durch
  $\sigma(x)(r) = rx$. Im kommutativen Diagramm
  \begin{equation*}
    \xymatrix{
      \operatorname{Hom}_R(R,M)\ar[r] & \operatorname{Hom}_R(I,M)\ar[r] &
      \operatorname{Ext}_R^1(R/I,M)\ar[r] & 0 \\
      M\ar[u]^{\cong}\ar[r]^(0.33){\sigma} &
      \operatorname{Hom}_R(I,IM)\ar@{^{(}->}[u]_{\iota_{\bullet}} 
    }
  \end{equation*}
  ist dann die obere Zeile exakt, sodass mit (\ref{1.2}, iii) die Behauptung
  folgt.
\end{Beweis}

\begin{Bemerkung}\label{1.7}
  Keine der beiden Bedingungen in (ii) kann man weglassen: Ist $R$ ein
  diskreter Bewertungsring und $0 \neq I \subsetneqq R$, so ist $\sigma$
  surjektiv, aber $\operatorname{Ext}_R^1(R/I,R) \neq 0$. Ist $R$ ein
  1-dimensionaler CM-Ring, aber kein diskreter Bewertungsring, so ist $R$
  $\mathfrak{m}$-ausgezeichnet, aber
  $\operatorname{Ext}_R^1(R/\mathfrak{m},R) \neq 0$.\\
  Ist aber $R$ beliebig und $I = Rc$ zyklisch, ist $\sigma$ automatisch
  surjektiv: Zu $f \colon I \to IM$ gibt es ein $x\in M$ mit $f(c) = cx$,
  und damit folgt $\sigma(x) = f$. Dasselbe $c$ liefert einen Isomorphismus
  $\operatorname{Hom}_R(I,M) \xrightarrow{\cong} M[\Ann_R(I)]$, und daraus
  folgt $\operatorname{Ext}_R^1(R/I,M) \cong
  \dfrac{M[\Ann_R(I)]}{IM}$.\\\mbox{}\hfill$\Box$
\end{Bemerkung}

Wir wollen eine dritte Beschreibung von "`$M$ ist $I$-ausgezeichnet"' durch
geeignete Erweiterungen $M \subset X$ formulieren, genauer durch die
Bestimmung von $\gamma_I(M)$. Falls $X$ injektiv ist, geschieht das in
(\cite{018}~p.~6), falls $M = R$ ein Integrit\"atsring, $X = K$ der
Quotientenk\"orper von $R$ und $I \neq 0$ ist, geschieht das in einer Reihe
von Arbeiten (\cite{001}, \cite{003}, \cite{005}, \cite{006}).

\begin{Lemma}\label{1.8}
  Sei $Y$ ein Untermodul von $X$ und $\alpha \colon Y :_X I \to \Hom_R(I,Y)$
  definiert durch $\alpha(u)(r) = ru$. Dann gilt:
  \begin{myenum}
    \item $\Hom_R(R/I, X) = 0 \implies \alpha$ injektiv.
    \item $\operatorname{Ext}_R^1(R/I,X) = 0 \implies \alpha$ surjektiv.
  \end{myenum}
\end{Lemma}

\begin{Beweis}
  F\"ur jedes $u \in X$ sei $h_u \colon R \to X$, $r \mapsto ru$, sodass
  $\alpha(u)$ die Einschr\"ankung von $h_u$ ist. Bei (1) ist
  $\operatorname{Ke} \alpha = \{u \in Y :_X I \,|\, ru = 0 \text{ f\"ur alle }
  r \in I\} = (Y :_X I)[I]$. Bei (2) ist $\Hom_R(R,X) \to \Hom_R(I,X)$
  surjektiv, sodass man zu jedem $f \in \Hom_R(I,Y)$ ein kommutatives Quadrat
  \begin{equation*}
    \xymatrix{
      I\ar[d]_f \ar@{}|{\textstyle\subset}[r] &  R\ar@{-->}[d]^{h_u} \\
      Y \ar@{}|{\textstyle\subset}[r] &  X
    }
  \end{equation*}
  erh\"alt, und aus $f(r) = ru$ f\"ur alle $r \in I$ folgt $u \in Y :_X I$,
  $\alpha(u) = f$.
\end{Beweis}

\begin{Proposition}\label{1.9}
  Sei $M \subset X$ eine Erweiterung mit $\operatorname{Ext}_R^1(R/I,X) =
  0$. Dann gilt:
  \begin{abc}
  \item $\gamma_I(M) = I(M :_X I)$.
  \item $M$ $I$-ausgezeichnet $\iff$ $(IM) :_X I = M :_X I$.
  \end{abc}
\end{Proposition}

\begin{Beweis}
  Auch ohne Voraussetzung an $X$ hat man das kommutative Diagramm
  \begin{equation*}
    \xymatrix{
      M\ar[dr]_{\sigma}\ar@{}|(0.4){\textstyle\subset}[r] & (IM) :_X I\ar[d]^{\alpha_1}\ar@{}|{\textstyle\subset}[r] & M :_X I\ar[d]^{\alpha_2} \ar@{}|{\textstyle\subset}[r]& X \\
      & \Hom_R(I, IM)\ar@{^{(}->}[r]_{\iota_{\bullet}} & \Hom_R(I,M) &
    }
  \end{equation*}
  in dem $\alpha_1$ und $\alpha_2$ wie in (\ref{1.8}) definiert seien.\\
  Bei (a) ist "`$\supset$"' klar, denn $I(M :_X I)$ ist ein $I$-generierter
  Untermodul von $M$, also enthalten in $\gamma_I(M)$. F\"ur "`$\subset$"'
  verwenden wir die Voraussetzung $\operatorname{Ext}_R^1(R/I,X) = 0$: Zu
  jedem $f \colon I \to M$ gibt es, weil $\alpha_2$ surjektiv ist, ein $u
  \in M :_X I$ mit $f = \alpha_2(u)$, also $f(r) = ru \in I(M :_X I)$ f\"ur
  alle $r \in I$, $\operatorname{Bi} f \subset I(M :_X I)$.\\
  Bei (b) gilt "`$\Rightarrow$"' wieder ohne Voraussetzung an $X$, denn aus
  $u \in M :_X I$ folgt
  \begin{equation*}
    \xymatrix{
      I\ar@{-->}[d]_f
      \ar@{}|{\textstyle\subset}[r]\ar@{}|{\textstyle{\circlearrowleft}}[dr] &
      R\ar@{-->}[d]^{h_u} \\ 
      M \ar@{}|{\textstyle\subset}[r] &  X ,     
    }
  \end{equation*}
  also nach (\ref{1.2}) $f(r) = ru \in IM$ f\"ur alle $r \in I$, d.\,h. $u \in
  (IM) :_X I$. Mit $\operatorname{Ext}_R^1(R/I,X) = 0$ ist auch
  "`$\Leftarrow$"' klar, denn dann ist im Diagramm $\alpha_2$, also auch
  $\iota_{\bullet}$ surjektiv, d.\,h. (\ref{1.2}, iii) erf\"ullt.
\end{Beweis}

\begin{Bemerkung}\label{1.10}
  Gilt sogar $\Hom_R(R/I,X) = \operatorname{Ext}_R^1(R/I,X) = 0$, liefert die
  exakte Folge $0 \to  \Hom_R(R/I,X/M) \to \operatorname{Ext}_R^1(R/I,M) \to
  0$ sofort $\operatorname{Ext}_R^1(R/I,M) \cong \dfrac{M :_X
    I}{M}$.\hfill$\Box$
\end{Bemerkung}

F\"ur jeden $R$-Modul $M$ hat man die Implikationen
\begin{equation*}
  M \text{ injektiv } \Rightarrow M \text{ ausgezeichnet } \Rightarrow M
  \text{ teilbar},
\end{equation*}
wobei $M$ wie in der Einleitung \emph{ausgezeichnet} hei{\ss}e, wenn es
$I$-ausgezeichnet ist f\"ur \emph{alle} Ideale $I$ (es gen\"ugt nach (\ref{1.5})
nur die zyklischen Ideale zu testen). Bei keiner dieser Implikationen gilt
die Umkehrung: Im ersten Fall w\"ahle $M = K/R$, wobei $R$ ein Integrit\"atsring
der Dimension $\geq 2$ sei (siehe \cite{011}~Theorem~5), im zweiten $M = R$,
wobei $R$ artinsch und $\dim_k (\operatorname{So}(R)) \geq 2$ sei (siehe
\ref{1.13}, Punkt~(h)).

Im folgenden Satz und seiner Folgerung sammeln wir einige Eigenschaften, die
f\"ur injektive Moduln wohlbekannt sind, aber auch noch f\"ur ausgezeichnete
Moduln gelten:
\begin{Satz}
  Ist $M$ ein ausgezeichneter $R$-Modul, so gilt:
  \begin{abc}
  \item F\"ur jedes Ideal $\mathfrak{a} \subset R$ ist die absteigende Folge
    $
      M \supset \mathfrak{a}M \supset \mathfrak{a}^2 M \supset \ldots
    $ station\"ar.
  \item Zu jedem $r \in R$ gibt es ein $e \geq 1$ mit $M[r^e] + rM = M$.
  \item $\operatorname{So}(M) \neq 0 \implies M$ ist treu.
  \item $\operatorname{So}(R) \neq 0 \implies \operatorname{So}(R) \cdot M =
    \operatorname{So}(M)$.
  \end{abc}
\end{Satz}

\begin{Beweis}
  (a) Diese Minimalbedingung f\"ur die Potenzen eines Ideals wurde in
  \cite{014} n\"aher untersucht und speziell in (Lemma~4.4) gezeigt, dass sie
  nur f\"ur \emph{zyklische} Ideale getestet werden muss. F\"ur jedes $r \in R$
  ist die aufsteigende Folge $\Ann_R(r) \subset \Ann_R(r^2) \subset \Ann_R(r^3)
  \subset \ldots$ station\"ar, und aus $\Ann_R(r^n) = \Ann_R(r^{n+1})$ f\"ur ein
  $n \geq 1$ folgt nach (\ref{1.2},~iv) $r^n M = r^{n+1}M$.\\
  (b) Aus $r^e M = r^{e+1} M$ wie in (a) folgt f\"ur jedes $x \in M$, dass
  $r^e x = r^{e+1}y$ ist mit $y \in M$, also $x-ry \in M[r^e]$.\\
  (c) Mit $I = \Ann_R(M)$ ist $IM = 0$, also nach Voraussetzung $\gamma_I(M)
  = 0$. Wegen $\mathfrak{m} \in \Ass(M)$ folgt nach (\ref{1.1}) $I=0$.\\
  (d) Nur "`$\supset$"' ist zu zeigen: Mit irgendeinem $0 \neq r \in
  \operatorname{So}(R)$ folgt $\Ann_R(r) = \mathfrak{m}$, also
  $\operatorname{So}(M) = M[\mathfrak{m}] = M[\Ann_R(r)] = rM$ nach
  (\ref{1.2}, iv) und damit die Behauptung. 
\end{Beweis}

\begin{Bemerkung}\label{1.12}
  Ein $R$-Modul $M \neq 0$ besitzt nach (\cite{014} Satz 4.2) genau dann eine
  Kompositionsreihe mit \emph{koprim\"aren} Faktoren, wenn $M$ die Bedingung
  (a) erf\"ullt und $\operatorname{Koass(M)}$ endlich ist.\hfill$\Box$
\end{Bemerkung}

Bei einem beliebigen $R$-Modul $M$ bezeichnen wir den gr\"o{\ss}ten radikalvollen
Untermodul mit $P(M)$, und $M/P(M)$ hei{\ss}t dann der \emph{reduzierte} Anteil
von $M$. War jetzt $M$ ausgezeichnet, folgt mit (a) bereits $P(M) =
\mathfrak{m}^n M$ f\"ur ein $n \geq 1$. (Selbst wenn $M$ injektiv war, muss
weder $P(M)$ noch $M/P(M)$ ausgezeichnet sein, z.\,B. wenn $M = E$ ist und
$H_{\mathfrak{m}}^0(R)$ einfach $\neq R$.)

\begin{Folgerung}\label{1.13}
  Ist $M$ ein \emph{reduzierter} ausgezeichneter $R$-Modul, so gilt:
  \begin{abc}
  \setcounter{abccount}{4}  
  \item $M$ ist schwach-flach.
  \item $\operatorname{So}(R) \cdot M$ ist gro{\ss} in $M$.
  \item $M \neq 0 \implies M$ ist treu und $R$ artinsch.
  \item $M$ zyklisch $\neq 0 \implies M \cong R$ und $R$ ist ein QF-Ring.
  \end{abc}
\end{Folgerung}

\begin{Beweis}
  (g) Nach der Vorbemerkung ist $\mathfrak{m}^n M = 0$, d.\,h.
  $\mathfrak{m}^n \subset \Ann_R(M)$ f\"ur ein $n \geq 1$. Falls $M \neq 0$,
  d.\,h. $\operatorname{So}(M) \neq 0$ ist, folgt nach (c) $M$ treu,
  $\mathfrak{m}^n = 0$, $R$ artinsch.\\
  (h) Nach eben ist $M \cong R$ und $R$ artinsch. Bleibt zu zeigen, dass
  $\operatorname{So}(R)$ einfach ist: Mit irgendeinem $0 \neq r \in
  \operatorname{So}(R)$ gilt $\Ann_R(r) = \mathfrak{m}$ und wie in (d) folgt
  $\operatorname{So}(R) = (r)$.\\
  (f) Nach eben ist $M$ halbartinsch, d.\,h. $\operatorname{So}(M)$ gro{\ss} in
  $M$. Mit (d) folgt die Behauptung.\\
  (e) Ein $R$-Modul $M$ hei{\ss}t nach (\cite{016}) schwach-flach, wenn f\"ur
  jeden Epimorphismus $g \colon Y \twoheadrightarrow M$ gilt, dass
  $\operatorname{Ke} g$ abgeschlossen in $Y$ ist. In (Satz~1.1) wird
  gezeigt, dass das \"aquivalent ist mit $\Ann_R(\mathfrak{p}) \cdot M$ gro{\ss}
  in $M[\mathfrak{p}]$ f\"ur alle $\mathfrak{p} \in \Ass(M)$. Weil in unserem
  Fall $M$ halbartinsch ist, folgt mit (f) die Behauptung.
\end{Beweis}

Der Spezialfall $M=R$ wird in (\cite{009} Remark~3.2) angek\"undigt:
\begin{Folgerung}\label{1.14} \emph{(Lindo)}
  Genau dann ist $R$ ausgezeichnet, wenn $R$ ein QF-Ring ist.
\end{Folgerung}

\begin{Bemerkung}\label{1.15}
  Aus Punkt (f) folgt f\"ur jeden reduzierten ausgezeichneten $R$-Modul $M$
  etwas allgemeiner: Ist $M$ Untermodul von $\mathfrak{m}X$, folgt $M = 0$.
  Zusammen mit (h) erh\"alt man also f\"ur jedes Ideal $0 \neq \mathfrak{a}
  \subsetneqq R$, dass weder $\mathfrak{a}$ noch $R/\mathfrak{a}$
  ausgezeichnet sein kann.\hfill$\Box$
\end{Bemerkung}

\section{Der Spezialfall $M =  R$}

Ein Ideal $I \subset R$ hei{\ss}e \emph{gut} (strong in \cite{001}, \cite{005},
trace ideal in \cite{008}), wenn $R$ $I$-ausgezeichnet ist, d.\,h. $I =
\gamma_I(R)$. Jede Summe von guten Idealen ist wieder gut (\ref{1.5}), und
aus $I \subset \gamma_I(R) \subset \Ann_R \Ann_R(I)$ folgt, dass jedes
Annullatorideal gut ist.

\begin{Lemma}\label{2.1}
  \begin{abc}
  \item Stets ist $\gamma_I(R)$ gut.
  \item $I \cong J \implies \gamma_I(R) = \gamma_J(R)$.
  \item $I$ gut und $\mathfrak{a}$ ein beliebiges Ideal $\implies I :_R
    \mathfrak{a}$ gut.
  \end{abc}
\end{Lemma}

\begin{Beweis}
  (a) Das wird in (\cite{008} Proposition~2.8, iv) gezeigt: $J =
  \gamma_I(R)$ ist $I$-generiert, sodass auch $\gamma_J(R)$ $I$-generiert
  ist, und in $\gamma_J(R) \subset J$ gilt nat\"urlich Gleichheit.\\
  (b) Allein aus $I \twoheadrightarrow J$ folgt f\"ur jeden $R$-Modul $M$,
  dass $\gamma_J(M)$ $I$-generiert ist, also $\gamma_J(M) \subset
  \gamma_I(M)$.\\
  (c) Zeigen wir f\"ur $J = I :_R \mathfrak{a}$ Punkt (ii) in (\ref{1.2}): Ist
  $f \colon J \to R$, folgt f\"ur alle $x \in J$, dass $\mathfrak{a} x \subset
  I$, also nach Voraussetzung $f(\mathfrak{a}x) \subset I$ ist,
  $\mathfrak{a}f(x) \subset I$, $f(x) \in  J$.
\end{Beweis}

\begin{Folgerung}
  Ist $I \cong J$, $I$ gut und artinsch, folgt bereits $I = J$.
\end{Folgerung}

\begin{Beweis}
  Aus $\text{L\"ange}\,(J) = \text{L\"ange}\,(\gamma_J(R)) < \infty$
  folgt, dass $J$ gut ist, also $J = I$. (Ohne "`artinsch"' muss $J$
  nicht einmal gut sein, siehe das \"ubern\"achste Beispiel.)
\end{Beweis}

\medskip

\noindent
\textbf{Beispiele} f\"ur $\gamma_I(R)$: (1) Ist $0 \neq I \subset
\operatorname{So}(R)$, so folgt $\gamma_I(R) = \operatorname{So}(R)$ (denn
aus $I \twoheadrightarrow R/\mathfrak{m}$ folgt $\operatorname{So}(R)
\subset \gamma_I(R)$, und "`$\supset$"' ist klar). (2) Ist das Ideal $I$
zyklisch und regul\"ar, so folgt $\gamma_I(R) = R$ (denn es ist $I \cong R$).
(3) Sei $R$ ein Integrit\"atsring und $I = \mathfrak{p} \cdot Rx$, wobei
$\mathfrak{p}$ ein Primideal der H\"ohe $n \geq 2$ und $0 \neq x \in
\mathfrak{p}$ sei. Dann ist $h(I) = 1$ und $h(\gamma_I(R)) = n$ (denn nach
(\cite{007}~p.~43, Ex.~26) ist $\mathfrak{p}$ gut, also $\mathfrak{p} =
\gamma_{\mathfrak{p}}(R) = \gamma_I(R)$).

\begin{Lemma}\label{2.3}
  Sei $M$ ein $R$-Modul und $\mathfrak{p} \in \operatorname{Spec}(R)$.
  \begin{abc}
  \item Ist $M$ $I$-ausgezeichnet, so ist auch $M_{\mathfrak{p}}$
    $IR_{\mathfrak{p}}$-ausgezeichnet.
  \item Ist $M_{\mathfrak{p}}$ $IR_{\mathfrak{p}}$-ausgezeichnet und
    operieren alle Elemente von $S = R \setminus \mathfrak{p}$ auf
    $M/IM$ injektiv, so ist auch $M$ $I$-ausgezeichnet.
  \end{abc}
\end{Lemma}

\begin{Beweis}
  (a) F\"ur jeden $R_{\mathfrak{p}}$-Homomorphismus $\alpha \colon
  IR_{\mathfrak{p}} \to M_{\mathfrak{p}}$ m\"ussen wir nach (\ref{1.2}, ii)
  $\operatorname{Bi} \alpha \subset IR_{\mathfrak{p}} \cdot
  M_{\mathfrak{p}}$ zeigen. Mit $\omega \colon I_{\mathfrak{p}}
  \xrightarrow{\ \cong\ } IR_{\mathfrak{p}}$ ist $\alpha \circ \omega$ nach
  (\cite{013} Theorem~3.84) von der Form $f_{\mathfrak{p}}$ mit $f \colon I
  \to M$, und f\"ur jedes $r \in I$ hat man nach Voraussetzung $f(r) = \sum
  r_ix_i$ mit $r_i \in I$, $x_i \in M$, also $\alpha(\frac{r}{s}) =
  \frac{f(r)}{s} = \sum \frac{r_i}{1} \cdot \frac{x_i}{s} \in
  I R_{\mathfrak{p}} \cdot M_{\mathfrak{p}}$ wie verlangt.\\
  (b) F\"ur jeden $R$-Homomorphismus $f \colon I \to M$ m\"ussen wir
  $\operatorname{Bi} f \subset IM$ zeigen. Ist $\omega$ wie in (a), gilt f\"ur
  den $R_{\mathfrak{p}}$-Homomorphismus $\alpha = f_{\mathfrak{p}} \circ
  \omega^{-1}$ nach Voraussetzung $\operatorname{Bi} \alpha \subset
  IR_{\mathfrak{p}} \cdot M_{\mathfrak{p}}$. Weil also $\alpha(\frac{r}{s}) =
  \frac{f(r)}{s}$ von der Form $\sum \frac{r_i x_i}{s'}$ ist, alle $r_i \in
  I$, $x_i \in M$, folgt $s''(s' f(r) - \sum sr_ix_i) = 0$, $s'' s' f(r)
  \in IM$, und weil $S$ auf $M/IM$ injektiv operiert, endlich $f(r) \in IM$.
\end{Beweis}

Weil $R$ genau dann ein diskreter Bewertungsring ist, wenn $\mathfrak{m}
\cong R$, also $\gamma_{\mathfrak{m}}(R) = R$, also $\mathfrak{m}$
\emph{nicht} gut ist, folgt mit (\ref{2.3}) sofort:
\begin{Satz}\label{2.4}
  Ein Primideal $\mathfrak{p}$ ist genau dann gut, wenn $R_{\mathfrak{p}}$
  \emph{kein} diskreter Bewertungsring ist.
\end{Satz}

Die in (\ref{1.6}) definierte Abbildung $\sigma \colon M \to \Hom_R(I,IM)$
wird jetzt zu einem Ringhomomorphismus $R \to \operatorname{End}_R(I)$, bei
dem das Ziel i. Allg. \emph{nicht} kommutativ ist. Z.\,B. liefert jede
Zerlegung $I = I_1 \oplus I_2$ mit $\Hom_R(I_1,I_2) \neq 0$ zwei $f,g \in
\operatorname{End}_R(I)$ mit $g \circ f \neq f \circ g$ (siehe \cite{008}
Remark~3.4) und damit folgt
\begin{Lemma}\label{2.5}
  Genau dann ist $R$ ein QF-Ring, wenn $\operatorname{So}(R) \neq 0$ ist und
  alle $\operatorname{End}_R(I)$ kommutativ sind.
\end{Lemma}

\begin{Beweis}
  "`$\Rightarrow$"' ist klar, denn $R$ ist artinsch und nach (\ref{1.6})
  sind alle $\sigma \colon R \to \operatorname{End}_R(I)$ sogar surjektiv.\\
  "`$\Leftarrow$"': Mit einem einfachen Untermodul $U \subset R$ und einem
  maximalen Element $V_0$ in der Menge $\{V \subset R \,|\, V \cap U = 0\}$
  ist $\operatorname{End}_R(V_0 \oplus U)$ kommutativ, also nach der
  Voraussetzung $V_0 = 0$, $U$ gro{\ss} in $R$, $R$ artinsch,
  $\operatorname{So}(R)$ einfach.
\end{Beweis}

\begin{Bemerkung}\label{2.6}
  Die folgenden Aussagen sind, falls $R$ ein Integrit\"atsring mit
  Quotientenk\"orper $K$ und $I \neq 0$ ist, wohlbekannt. Sie gelten
  allgemeiner f\"ur jeden Ring $R$, wenn $K$ der totale Quotientenring von $R$
  ist und $I$ ein \emph{regul\"ares} Ideal. Dann ist n\"amlich
  $\operatorname{Ext}_R^i(R/I,K) = 0$ f\"ur alle $i \geq 0$, und die
  Ergebnisse von (\ref{1.8}) bis (\ref{1.10}) zeigen:\\[2ex]
  Mit $I^{-1} = R :_K I$ sind im kommutativen Diagramm
  \begin{equation*}
    \xymatrix{
      R\ar[dr]_{\sigma}\ar@{}|(0.4){\textstyle\subset}[r] & I :_K I\ar[d]^{\alpha_1}\ar@{}|{\textstyle\subset}[r] & I^{-1}\ar[d]^{\alpha_2} \ar@{}|{\textstyle\subset}[r]& K \\
      & \Hom_R(I, I)\ar[r] & \Hom_R(I,R) &
    }
  \end{equation*}
  $\alpha_1$ und $\alpha_2$ Isomorphismen, $\gamma_I(R) = I \cdot I^{-1}$,
  $I$ gut $ \Leftrightarrow I :_K I = I^{-1}$ und
  $\operatorname{Ext}_R^1(R/I,R) \cong I^{-1}/R$.\\
  Weil $I$ regul\"ar, also $R \subset I :_K I$ eine endliche Ringerweiterung
  ist, gilt insbesondere: Ist $R$ in seinem totalen Quotientenring $K$ ganz
  abgeschlossen, folgt f\"ur jedes regul\"are gute Ideal $I$ bereits
  $\operatorname{Ext}_R^1(R/I,R) = 0$.\hfill$\Box$
\end{Bemerkung}

Sei nun $R$ ein 1-dimensionaler CM-Ring. Dann lassen sich einige Ergebnisse
von Kapitel~XII und XIII in \cite{012} so zusammenfassen: Ist $e$ die
Multiplizit\"at von $R$ und $v(\mathfrak{a}) = \dim_k(\mathfrak{a} /
\mathfrak{m}\mathfrak{a})$ die minimale Erzeugendenanzahl eines Ideals
$\mathfrak{a}$, so gibt es eine nat\"urliche Zahl $\nu \geq 0$ mit
$v(\mathfrak{m}^n) = e$ f\"ur alle $n \geq \nu$, $v(\mathfrak{m}^n) < e$ f\"ur
alle $n < \nu$ (12.10). F\"ur den ersten Nachbarschaftsring $\Lambda =
\{\frac{b}{a} \in K \,|\, b \in \mathfrak{m}^s\ (s \geq 1)$ und $a$ ist ein
superfizielles Element vom $\operatorname{Grad} s\}$ gilt dann
$\mathfrak{m}^n \cdot (\mathfrak{m}^n)^{-1} = \Lambda^{-1}$ f\"ur alle $n \geq
\nu$ (12.12). Falls $v(\mathfrak{m}) = 2$ ist, gilt sogar $\nu = e-1$ und
$\Lambda = (\mathfrak{m}^{e-1})^{-1}$, $\Lambda^{-1} = \mathfrak{m}^{e-1}$
(13.8). Damit folgt
\begin{Satz}\label{2.7} \emph{(Matlis)}
  Sei $R$ ein 1-dimensionaler CM-Ring und seien $\nu$, $\Lambda$ wie in der
  Vorbereitung. F\"ur alle Ideale $I \in
  \{\mathfrak{m}^{\nu},\mathfrak{m}^{\nu+1},\mathfrak{m}^{\nu+2},\dotsc\}$
  gilt dann \[\gamma_I(R) = \Lambda^{-1}.\] Falls $v(\mathfrak{m}) = 2$ ist,
  gilt weiter $\Lambda^{-1} = \mathfrak{m}^{\nu}$.
\end{Satz}

\section{Dualisierung}

Ist $\kappa_I(M)$ das kleinste Element in der Menge $\{V \subset M \,|\, M/V
\text{ ist } I^\circ\text{-kogeneriert} \}$, d.\,h. $\kappa_I(M) = \bigcap
\{\operatorname{Ke} f \,|\, f \in \Hom_R(M, I^\circ)\}$, so gilt wieder nach (\cite{018}~p.~6)
\begin{equation*}
  \Ann_R(I) \cdot M \subset \kappa_I(M) \subset M[I],
\end{equation*}
und wir sagen, $M$ sei \emph{$I$-koausgezeichnet}, wenn $\kappa_I(M) = M[I]$
ist, \emph{koausgezeichnet}, wenn diese Gleichung f\"ur alle Ideale $I$ von
$R$ gilt. Jeder flache $R$-Modul ist koausgezeichnet (hier stimmen alle drei
Untermoduln \"uberein), und falls $I$ zyklisch ist, gilt $\Ann_R(I) \cdot M =
\kappa_I(M)$ (\cite{018}~p.~5,~6). Die n\"achsten zwei Aussagen (\ref{3.1})
und (\ref{3.2}) ergeben sich, durch das Matlis-Duale $M^\circ = \Hom_R(M,
E)$ und die Isomorphismen $\operatorname{Ext}_R^i (M,N^\circ) \cong
\operatorname{Tor}_i^R(M,N)^\circ$ (\cite{004}~p.~120, Proposition~5.1),
sofort aus Teil~(1), sodass wir ihren Beweis dem Leser \"uberlassen (siehe
auch 3.7,~a).

\begin{Lemma}\label{3.1}
  F\"ur einen $R$-Modul $M$ und ein Ideal $I \subset R$ sind \"aquivalent:
  \begin{iii}
    \item $\kappa_I(M) = M$.
    \item $M \otimes_R I = 0$.
    \item $I_{\mathfrak{p}} = 0$ f\"ur alle $\mathfrak{p} \in
      \operatorname{Koass}(M)$.
  \end{iii}
\end{Lemma}

\begin{Proposition}\label{3.2}
  F\"ur einen $R$-Modul $M$ und ein Ideal $I \subset R$ sind \"aquivalent:
  \begin{iii}
    \item $M$ ist $I$-koausgezeichnet, d.\,h. $\kappa_I(M) = M[I]$.
    \item F\"ur jeden Homomorphismus $f \colon M \to I^\circ$ gilt $M[I]
      \subset \operatorname{Ke} f$, d.\,h.
      \begin{equation*}
        \xymatrix{
          M \ar_f[d]\ar[r]\ar@{}|(.3){\textstyle{\circlearrowleft}}[dr] &
          M/M[I]\ar@{-->}[dl] \\
          I^\circ &
        }
      \end{equation*}
    \item  $M[I] \otimes_R I \xrightarrow{\ \iota_\bullet\ } M \otimes_R I$
      ist die Nullabbildung.
    \end{iii}
    War das Ideal $I$ zyklisch, ist das weiter \"aquivalent mit
    \begin{iii}
      \setcounter{iiicount}{3}
      \item $\Ann_R(I) \cdot M = M[I]$.
    \end{iii}
\end{Proposition}

F\"ur den (\ref{1.6}) entsprechenden Zusammenhang zwischen "`$M$ ist
$I$-koausgezeichnet"' und "`$\operatorname{Tor}_1^R(M, R/I) =0$"' definieren
wir die kanonische Abbildung $\beta$ durch
\begin{equation*}
  \xymatrix{
    M[I] \otimes_R I\ar[r]\ar_0[dr] & M \otimes_R
    I\ar[r]\ar[d]\ar@{}|(.33){\textstyle{\circlearrowleft}}[dr] &
    \dfrac{M}{M[I]} \otimes_R I \ar[r]\ar@{-->}[dl]^{\beta} & 0, \\
    & M & & &
  }
\end{equation*}
also $\beta(\overline{x} \otimes r) = rx$ f\"ur $x \in M$, $r \in I$. Im
kommutativen Diagramm
\begin{equation*}
  \xymatrix{
    0 \ar[r] & \operatorname{Tor}_1^R(M,R/I) \ar[r] & M \otimes_R
    I\ar@{->>}^(0.4){\nu_\bullet}[d] 
    \ar^{\iota_{\bullet}}[r] & M \otimes_R R \ar^{\cong}[d] \\
    & & \dfrac{M}{M[I]} \otimes_R I \ar_{\beta}[r] & M
  }
\end{equation*}
ist die obere Zeile exakt, sodass mit (\ref{3.2},~iii) sofort folgt:
\begin{Proposition}\label{3.3}
  F\"ur einen $R$-Modul $M$ und ein Ideal $I \subset R$ sind \"aquivalent:
  \begin{iii}
    \item $\operatorname{Tor}_1^R(M, R/I) = 0$.
    \item $M$ ist $I$-koausgezeichnet und die kanonische Abbildung $\beta
      \colon \dfrac{M}{M[I]} \otimes_R I \to M$ injektiv.
  \end{iii}
\end{Proposition}

\begin{Bemerkung}\label{3.4}
  War das Ideal $I = Rc$ zyklisch, ist $\beta$ automatisch injektiv: Jedes
  Element von $\dfrac{M}{M[I]} \otimes_R I$ ist dann von der Form
  $\overline{x} \otimes c$, sodass aus $\beta(\overline{x} \otimes c) = 0$
  folgt $cx = 0$, $x \in M[I]$, $\overline{x} = 0$. Mit demselben $c$ erh\"alt
  man einen Isomorphismus $M \otimes_R I \xrightarrow{\ \cong\ } M/\Ann_R(I)
  \cdot M$ und damit $\operatorname{Tor}_1^R(M,R/I) \cong
  \dfrac{M[I]}{\Ann_R(I) \cdot M}$.\hfill$\Box$
\end{Bemerkung}

Wir wollen die drei Bedingungen in der letzten Proposition durch ein
Beispiel illustrieren:

\medskip

\noindent
\textbf{Beispiel} Sei $K$ der totale Quotientenring von $R$ und $I \subset
R$ regul\"ar. Dann gilt:
\begin{myenum}
  \item $\operatorname{Tor}_1^R(K/R, R/I) = 0 \iff I = R$.
  \item $K/R$ $I$-koausgezeichnet $\iff I :_K I = I^{-1}$.
  \item $\beta \colon \dfrac{K/R}{(K/R)[I]} \otimes_R I \to K/R$ injektiv
    $\iff I \cong R$.
\end{myenum}

\begin{Beweis}
  (1) F\"ur jeden $R$-Modul $N$ ist $\operatorname{Tor}_1^R(K/R,N) \cong
  T(N)$, sodass die linke Seite \"aquivalent ist mit $T(R/I) = 0$, d.\,h.
  $I=R$.\\
  (2) Weil ${}_RK$ flach ist, gilt $\kappa_I(K/R) = (I :_K I)/R$ nach
  (\cite{018}~p.~6), au{\ss}erdem $(K/R)[I] = (R :_K I)/R = I^{-1}/R$, sodass
  die linke Seite \"aquivalent ist mit $I :_K I = I^{-1}$.\\
  (3) Der Isomorphismus $\omega\colon \dfrac{K/R}{(K/R)[I]} \otimes_R I
  \xrightarrow{\ \cong\ } \dfrac{K}{I^{-1}} \otimes_R I \xrightarrow{\
    \cong\ } K/I \cdot I^{-1}$ induziert den zu $\beta$ \"aquivalenten
  Homomorphismus $\beta' = \beta \circ \omega^{-1} \colon K/I \cdot I^{-1}
  \to K/R$, und weil das gerade die kanonische Abbildung ist, folgt
  $\operatorname{Ke} \beta \cong \operatorname{Ke} \beta' = R/I \cdot
  I^{-1}$. Die linke Seite von (3) ist also \"aquivalent mit $I \cdot I^{-1} =
  R$, $R$ $I$-generiert, $I \cong R$. 
\end{Beweis}

Die zu (\ref{1.9}) duale Aussage lautet:
\begin{Proposition}\label{3.5}
  Sei $B \subset A$ eine Erweiterung mit $\operatorname{Tor}_1^R(A,R/I) =
  0$. Dann gilt:
  \begin{abc}
  \item $\kappa_I(A/B) = \bigl((IB) :_A I\bigr)/B$.
  \item $A/B$ $I$-koausgezeichnet $\iff (IB) :_A I = B :_A I$.
  \end{abc}
\end{Proposition}

\begin{Beweis}
  (b) ist ein Spezialfall von (a), und mit $B^\ast = (IB) :_A I$ ist
  $\frac{A}{IB}/ \frac{A}{IB}[I] \cong \frac{A}{B}/\frac{B^\ast}{B}$ nach
  (\cite{017}~Prop.~3.1) $I^\circ$-kogeneriert, also stets $\kappa_I(A/B)
  \subset B^\ast/B$. Zur Umkehrung m\"ussen wir f\"ur jedes $f \colon A/B \to
  I^\circ$ zeigen, dass $B^\ast/B \subset \operatorname{Ke} f$ ist, mit
  $\overline{f} \colon I \to (A/B)^\circ$ also
  $\operatorname{Bi}\overline{f} \subset \Ann_{(A/B)^\circ}(B^\ast/B)$, mit
  $\nu\colon A \to A/B$ also $\nu^\circ(\operatorname{Bi}\overline{f})
  \subset \Ann_{A^\circ}(B^\ast)$.\\
  Mit $X = A^\circ$ und $M = \Ann_{A^\circ}(B)$ folgt nun, weil nach
  Voraussetzung $\operatorname{Ext}_R^1(R/I,X) \cong
  \operatorname{Tor}_1^R(A,R/I)^\circ = 0$ ist, $\gamma_I(M) = I(M :_X I)$
  nach (\ref{1.9}), $I(M :_X I) = \Ann_{A^\circ}(B^\ast)$ nach (\cite{017}
  Bemerkung~4.4), mit $\nu^\circ(\operatorname{Bi}\overline{f}) \subset
  \gamma_I(M)$ endlich die Behauptung.
\end{Beweis}

Weitere Parallelit\"aten zum Teil~(1) ergeben sich aus dem n\"achsten Satz, den
wir in (\cite{018} Folgerung~3.6) nur f\"ur einreihige Moduln zeigten.

\begin{Satz}\label{3.6}
  F\"ur jeden $R$-Modul $M$ gilt:
  \begin{abc}
  \item $\gamma_I(M^\circ) = \Ann_{M^\circ}(\kappa_I(M))$.
  \item $\kappa_I(M^\circ) = \Ann_{M^\circ}(\gamma_I(M))$.
  \end{abc}
\end{Satz}

\begin{Beweis}
  (a) Sei $V = \kappa_I(M)$. Bei "`$\supset$"' ist $M/V$
  $I^\circ$-kogeneriert, also $(M/V)^\circ \cong \Ann_{M^\circ}(V)$ nach
  (\cite{017} Proposition~3.1) $I$-generiert, und daraus folgt die
  Behauptung. Bei "`$\subset$"' m\"ussen wir f\"ur jeden Homomorphismus $g
  \colon I \to M^\circ$ zeigen, dass $\operatorname{Bi} g \subset
  \Ann_{M^\circ}(V)$ ist, mit $f \colon M \to I^\circ$ und $g =
  \overline{f}$ also $V \subset \operatorname{Ke} f$, und das ist klar.\\
  (b) Sei $M \subset X$ eine injektive Erweiterung und $A = X^\circ$, $B =
  \Ann_{X^\circ}(M)$. Dann ist die induzierte Abbildung $\omega \colon A/B
  \to M^\circ$ ein Isomorphismus und wir wollen zuerst $\kappa_I(A/B)$
  berechnen: Weil $A$ flach ist, gilt $\kappa_I(A/B) = B^\ast/B$ mit $B^\ast
  = (IB) :_A I$ nach (\ref{3.5}), weiter $B^\ast = \Ann_{X^\circ}(M_\ast)$
  mit $M_\ast = I(M :_X I)$ nach (\cite{017} Lemma~4.3), schlie{\ss}lich $M_\ast
  = \gamma_I(M)$ nach (\ref{1.9}), also zusammen
  \begin{equation*}
    \kappa_I(A/B) = \Ann_{X^\circ}(\gamma_I(M)) / \Ann_{X^\circ}(M).
  \end{equation*}
  Der Isomorphismus $\omega$ f\"uhrt dann $\kappa_I(A/B)$ in
  $\Ann_{M^\circ}(\gamma_I(M))$ \"uber.
\end{Beweis}

\begin{Folgerung}\label{3.7}
  \mbox{}
  \begin{abc}
  \item $M$ $I$-koausgezeichnet $\iff M^\circ$ $I$-ausgezeichnet.
  \item $M$ $I$-ausgezeichnet $\iff M^\circ$ $I$-koausgezeichnet.
  \end{abc}
\end{Folgerung}

\begin{Beweis}
  Klar mit den Gleichungen $\Ann_{M^\circ}(M[\mathfrak{a}]) =
  \mathfrak{a}\cdot M^\circ$ und $\Ann_{M^\circ}(\mathfrak{a} M) = M^\circ[\mathfrak{a}]$.
\end{Beweis}

\begin{Folgerung}\label{3.8}
  Ist $U$ ein reiner Untermodul von $M$, so gilt:
  \begin{abc}
  \item $\gamma_I(U) = U \cap \gamma_I(M)$.
  \item $\kappa_I(U) = U \cap \kappa_I(M)$.
  \end{abc}
\end{Folgerung}

\begin{Beweis}
  (a) Weil $\iota^\circ\colon M^\circ \to U^\circ$ ein zerfallender
  Epimorphismus ist, gilt das auch f\"ur $\kappa_I(M^\circ) \to
  \kappa_I(U^\circ)$, d.\,h. nach (\ref{3.6},~b) auch f\"ur
  $\bigl(M/\gamma_I(M)\bigr)^\circ \to \bigl(U/\gamma_I(U)\bigr)^\circ$.
  Also ist die urspr\"ungliche Abbildung $U/\gamma_I(U) \to M/\gamma_I(M)$
  eine Monomorphismus, d,\,h. $U \cap \gamma_I(M) \subset \gamma_I(U)$, und
  "`$\supset$"' ist klar.\\
  (b) Sei $x \in U \cap \kappa_I(M)$ und $g \colon U \to I^\circ$. Weil
  $I^\circ$ rein-injektiv ist, gibt es ein $f \colon M \to I^\circ$ mit
  $f|_U =g$, sodass aus $x \in \kappa_I(M)$ folgt $f(x) = 0$, aus $x \in U$
  schlie{\ss}lich $g(x) = 0$.
\end{Beweis}

\begin{Bemerkung}\label{3.9}
  F\"ur jeden reinen Untermodul $U$ von $M$ kann man entsprechend\\
  $\gamma_I(M/U) = (\gamma_I(M) + U)/U$ und $\kappa_I(M/U) = (\kappa_I(M) +
  U)/U$ zeigen, indem man f\"ur jeden reinen Epimorphismus $A
  \twoheadrightarrow B$ den zerfallenden Monomorphismus $B^\circ
  \hookrightarrow A^\circ$ betrachtet und wieder (\ref{3.6})
  verwendet.\hfill$\Box$
\end{Bemerkung}

\medskip

\noindent
\textbf{Beispiel} Sei $R$ vollst\"andig, $K$ der totale Quotientenring von $R$
und $I \subset R$ regul\"ar. Dann gilt:
\begin{myenum}
  \item $\operatorname{Tor}_1^R(E, R/I) = 0 \iff R = I^{-1}$.
  \item $E$ $I$-koausgezeichnet $\iff I :_K I = I^{-1}$.
  \item $\beta\colon \dfrac{E}{E[I]} \otimes_R I \to E$ injektiv $\iff R=I
    :_K I$.
\end{myenum}

\begin{Beweis}
  (1) Wegen $\operatorname{Ext}_R^1(R/I,R) \cong
  \operatorname{Tor}_1^R(E,R/I)^\circ$ ist die linke Seite nach (\ref{2.6})
  \"aquivalent mit $R = I^{-1}$.\\
  (2) Nach (\ref{3.7},~b) ist die linke Seite \"aquivalent mit $R$
  $I$-ausgezeichnet, also wieder nach (\ref{2.6}) \"aquivalent mit $I :_K I =
  I^{-1}$.\\
  (3) F\"ur jeden $R$-Modul $M$ hat man einen kanonischen Isomorphismus
  $\omega$ mit
  \begin{equation*}
  \xymatrix{
  M^\circ\ar_{\beta^\circ}[dr]\ar^(0.4){\sigma}[rr] &
  \ar@{}|(0.4){\textstyle{\circlearrowleft}}[d]&
  \Hom_R(I, I \cdot M^\circ) \\ & \quad\left(\frac{M}{M[I]}\otimes_R
    I\right)^\circ\ar_\omega[ur] &
  }
  \end{equation*}
  sodass $\beta$ genau dann injektiv ist, wenn $\sigma$ surjektiv ist.
  Speziell f\"ur $M=E$ ist also die linke Seite in (3) \"aquivalent mit der
  Surjektivit\"at von $R \xrightarrow{\ \sigma\ } \Hom_R(I,I)$, d.\,h. mit $R
  = I :_K I$.
\end{Beweis}

Wie in Teil (1) jeder ausgezeichnete $R$-Modul teilbar war, ist jetzt jeder
koausgezeichnete $R$-Modul torsionsfrei (und beide Male gilt, falls $R$ ein
Integrit\"atsring ist, die Umkehrung). Diese Tatsachen und die im Folgenden
gesammelten Eigenschaften von koausgezeichneten $R$-Moduln lassen sich jetzt
ohne M\"uhe durch das Matlis-Duale $M^\circ$ und (\ref{3.7},~a) auf Teil~(1)
zur\"uckf\"uhren (f\"ur die schwach-injektiven $R$-Moduln in Punkt~(e)
siehe~\cite{015}), sodass wir nur noch die Hauptergebnisse anschreiben.

\begin{Satz}\label{3.10}
  Ist $M$ ein koausgezeichneter $R$-Modul, so gilt:
  \begin{abc}
  \item F\"ur jedes Ideal $\mathfrak{a} \subset R$ ist die aufsteigende Folge
    $M[\mathfrak{a}] \subset M[\mathfrak{a}^2] \subset M[\mathfrak{a}^3]
    \subset \ldots$ station\"ar.
  \item Zu jedem $r \in R$ gibt es ein $e \geq 1$ mit $r^e M \cap M[r] = 0$.
  \item $\mathfrak{m} M \neq M \implies M$ ist treu.
  \item $\operatorname{So}(R) \neq 0 \implies \mathfrak{m}M =
    M[\operatorname{So}(R)]$.
 \end{abc}
\end{Satz}

\begin{Folgerung}\label{3.11}
  Ist $M$ ein \emph{halbartinscher} koausgezeichneter $R$-Modul, so gilt:
  \begin{abc}
  \setcounter{abccount}{4}
  \item $M$ ist schwach-injektiv.
  \item $M[\operatorname{So}(R)]$ ist klein in $M$.
  \item $M \neq 0 \implies M$ ist treu und $R$ artinsch.
  \item $M$ kozyklisch $\neq 0 \implies M \cong R$ und $R$ ist ein QF-Ring. 
  \end{abc}
\end{Folgerung}


\begin{thebibliography}{99}
\bibitem{001} V. Barucci: \emph{Strongly divisorial ideals and complete
    integral closure of an integral domain}: J. Algebra 99 (1986) 132--142
\bibitem{002}
  N. Bourbaki: \emph{Algèbre commutative}: Hermann. Paris (1967)
\bibitem{003} J. P. Brennan -- W. V. Vasconcelos: \emph{On the structure of
    closed ideals}: Math. Scand. 88 (2001) 3--16
\bibitem{004} H. Cartan -- S. Eilenberg: \emph{Homological algebra}:
  Princeton Univ. Press (1956)
\bibitem{005} M. Fontana -- J. A. Huckaba -- I. J. Papick: \emph{Domains
    satisfying the trace property}: J. Algebra 107 (1987) 169--182
\bibitem{006} J. A. Huckaba -- I. J. Papick: \emph{When the dual of an ideal
    is a ring}: manuscripta math. 37 (1982) 67--85
\bibitem{007} I. Kaplansky: \emph{Commutative rings}: Chicago Univ. Press
  (1974)
\bibitem{008} H. Lindo: \emph{Trace ideals and centers of endomorphism rings
    of modules over commutative rings}: J. Algebra 482 (2017) 102--130
\bibitem{009} H. Lindo: \emph{Self-injective commutative rings have no
    nontrivial rigid ideals}: arXiv 1710.01973 (2017) 1--6
\bibitem{010} H. Lindo -- N. Pande: \emph{Trace ideals and the Gorenstein
    property}: arXiv 1802.0649 (2018) 1--6
\bibitem{011} E. Matlis: \emph{Some properties of noetherian domains of
    dimension one}: Canadian J. Math. 13 (1961) 569--586
\bibitem{012} E. Matlis: \emph{1-dimensional Cohen-Macaulay rings}: Springer
  LNM 327 (1973)
\bibitem{013} J. J. Rotman: \emph{An introduction to homological algebra}:
  Academic Press. New York (1979)
\bibitem{014} H. Z"oschinger: \emph{Moduln mit Koprim\"arzerlegung}: Bayer.
  Akad. Wiss. 2 (1990) 5--25
\bibitem{015} H. Z"oschinger: \emph{Schwach-injektive Moduln}: Per. Math.
  Hung. 52 (2006) 105--128
\bibitem{016} H. Z"oschinger: \emph{Schwach-flache Moduln}: Commun. Algebra
  41 (2013) 4393--4407
\bibitem{017}  H. Z"oschinger: \emph{\"Uber die von einem Ideal $I \subset R$
    erzeugten $R$-Moduln}: arXiv 1604.02349 (2016) 1--9
\bibitem{018}  H. Z"oschinger: \emph{\"Uber die von einem Ideal $I \subset R$
    erzeugten $R$-Moduln II}: arXiv 1705.03353 (2017) 1--8
\end{thebibliography}
\end{document}